\documentclass[11pt]{article}
\usepackage{amsmath}
\usepackage{latexsym}
\usepackage{graphicx}
\makeatletter
\newfont{\footsc}{cmcsc10 at 8truept}
\newfont{\footbf}{cmbx10 at 8truept}
\newfont{\footrm}{cmr10 at 10truept}
\makeatother
\newtheorem{theorem}{\bf Theorem}

\newtheorem{lemma}{\bf Lemma}
\begin{document}

\title{On an Asymptotic Series of Ramanujan}

\author{Yaming Yu\\
\small Department of Statistics\\[-0.8ex]
\small University of California\\[-0.8ex] 
\small Irvine, CA 92697, USA\\[-0.8ex]
\small \texttt{yamingy@uci.edu}}

\date{}
\maketitle

\begin{abstract}
An asymptotic series in Ramanujan's second notebook (Entry 10, Chapter 3) is concerned with the behavior of the expected 
value of $\phi(X)$ for large $\lambda$ where $X$ is a Poisson random variable with mean $\lambda$ and $\phi$ is a function 
satisfying certain growth conditions.  We generalize this by studying the asymptotics of the expected value of $\phi(X)$ when 
the distribution of $X$ belongs to a suitable family indexed by a convolution parameter.  Examples include the binomial, 
negative binomial, and gamma families.  Some formulas associated with the negative binomial appear new. 


\end{abstract}

\section{An Asymptotic Series of Ramanujan}
A version (modified from \cite{B}) of Entry 10 in Chapter 3 of Ramanujan's second notebook reads

\begin{theorem}
\label{thm1}
Let $\phi(x),\ x\in [0,\infty),$ denote a function of at most polynomial growth as $x$ tends to $\infty$.  Suppose there 
exist constants $x_0>0$ and $A\geq 1$, and a function $G(x)$ of at most polynomial growth as $x\rightarrow \infty$ 
such that for each nonnegative integer $m$ and all $x>x_0$, the derivatives $\phi^{(m)}(x)$ exist and satisfy
\begin{equation}
\label{cond}
\left| \frac{\phi^{(m)}(x)}{m!}\right|\leq G(x)\left(\frac{A}{x}\right)^m.
\end{equation}
Assume that there exists a positive constant $c$ such that
\begin{equation}
\label{origin}
G(x)\gg e^{-c\sqrt{x}}
\end{equation}
as $x\rightarrow\infty$.  Put
$$\phi_\infty(x)=e^{-x} \sum_{k=0}^\infty \frac{x^k\phi(k)}{k!}.$$
Then for any fixed positive integer $M$,
\begin{equation}
\label{eqn1}
\phi_{\infty}(x)=\phi(x)+\sum_{n=2}^{2M-2} \sum_{k=\lfloor (n+1)/2\rfloor +1}^n b_{kn} x^{n-k+1} 
\frac{\phi^{(n)}(x)}{n!}+O\left(G(x)x^{-M}\right),
\end{equation}
as $x$ tends to $\infty$, where $\lfloor (n+1)/2\rfloor $ denotes the integer part of $(n+1)/2$ and the numbers $b_{kn}$ are 
defined recursively by
\begin{align*}
b_{kk} &=1, &k\geq 2;\\ 
b_{kn} &=0, &n<k\ {\rm or}\ n>2k-2;\\
b_{k+1, n+1}&=nb_{k,n-1}+(n-k+1)b_{kn}, &k\leq n\leq 2k-1.
\end{align*}
\end{theorem}

This result may seem hard to penetrate at first glance.  Its relevance, however, is easily appreciated through  interesting 
examples such as (\cite{B, B2})
\begin{equation}
\label{sqrt}
e^{-x}\sum_{k=0}^\infty \frac{\sqrt{k}x^k 
}{k!}=\sqrt{x}\left(1-\frac{1}{8x}-\frac{7}{128x^2}+O\left(\frac{1}{x^3}\right)\right)
\end{equation}
and 
$$ e^{-x} \sum_{k=0}^\infty \frac{x^k\log(k+1)}{k!}=\log(x)+\frac{1}{2x}+\frac{1}{12x^2}+O\left(\frac{1}{x^3}\right),$$
both valid as $x\rightarrow\infty$; by choosing $\phi(x)=\log\Gamma(x+1)$ in (\ref{eqn1}), an asymptotic formula for the Shannon entropy of the Poisson distribution can also be obtained (see \cite{EBBJ}).

The first goal of this note is a formal, probablistic derivation of Theorem \ref{thm1}.  The starting point is the observation 
that $\phi_\infty (x) = E\phi(U)$, where $U$ is a Poisson random variable with mean $x$ and $E$ denotes expectation.  Based on 
this we present in Section 2 a more general version (Theorem \ref{main}) of Theorem \ref{thm1} by considering $U$ distributed 
as some distribution other than the Poisson, e.g., a gamma distribution or a binomial distribution.  As illustrations, we 
derive asymptotic expansions for digamma functions and for inverse moments of certain positive random variables.  We prove 
Theorem \ref{main} in Section 3. 

Noting that $\phi_\infty (x) = E\phi(U)$ where $U$ has the ${\rm Po}(x)$ distribution, we expand $\phi(U)$ as a 
Taylor series
$$\phi(U)=\phi(x)+\sum_{n=1}^{2M-2} \frac{(U-x)^n\phi^{(n)}(x)}{n!}+\ldots,$$
and formally take the expectation term by term:
\begin{equation}
\label{ephi}
E\phi(U)=\phi(x)+\sum_{n=1}^{2M-2} \frac{E(U-x)^n\phi^{(n)}(x)}{n!}+\ldots
\end{equation}
The quantity $\mu_n=E(U-x)^n$ is the $n$th central moment of the ${\rm Po}(x)$ distribution.  The first few $\mu_n$'s are
$$(\mu_1, \mu_2, \mu_3, \mu_4)=(0, x, x, 3x^2+x),$$
and they obey the well-known recursion (see \cite{Z}, Lemma 3, for example)
$$\mu_{n+1}=x\left(\frac{{\rm d}\mu_n}{{\rm d}x}+n\mu_{n-1}\right),\quad n\geq 2,$$
from which we obtain, by comparing the coefficients of $x^{n-k+1}$ and using the definition of $b_{kn}$, 
\begin{equation}
\label{mu}
\mu_n=\sum_{k=0}^n b_{kn} x^{n-k+1},\quad n\geq 2.
\end{equation}
The double sum in (\ref{eqn1}) is the result of substituting (\ref{mu}) in (\ref{ephi}) and noting that $b_{kn}=0$ if $k\leq 
\lfloor (n+1)/2\rfloor$.  Based on this it is also clear that $\mu_n$ is a polynomial in $x$ of degree $\lfloor 
n/2\rfloor$, which implies that, given the condition (\ref{cond}), the ``leading term of the remainder'' in (\ref{ephi}),  
$$\frac{E(U-x)^{2M-1}\phi^{(2M-1)}(x)}{(2M-1)!}=\frac{\mu_{2M-1}\phi^{(2M-1)}(x)}{(2M-1)!},$$
is $O(G(x)x^{-M}).$

The above derivation is, of course, strictly formal.  However, it can be made rigorous under the stated conditions; see Berndt 
\cite{B} and Evans \cite{E}.  Berndt actually proved a modification of (\ref{eqn1}) where the order of summation over $k$ and 
$n$ on the right hand side is inverted and certain higher order terms of the resulting sum are absorbed in the $O(G(x) 
x^{-M})$ term. 

\section{A General Version}
The formal derivation in Section 1 suggests that it is possible to generalize Theorem \ref{thm1} if we let $U$ have a suitable 
distribution other than the Poisson.  Noting the key role played by the central moments of $U$, we give a version of Theorem 
\ref{thm1} by imposing conditions on the moment generating function (mgf) of $U$.  An introduction to moment generating 
functions can be found in probability texts such as Gut \cite{G}.  A useful property is that, if an mgf exists in a 
neighborhood of zero, then for all $m\geq 1$, the $m$th moment exists and can be obtained by differentiating the mgf $m$ 
times. 

\begin{theorem}
\label{main}
Let $\phi(x),\ x\in [0,\infty),$ denote a Borel measurable function that can be bounded in absolute value by a polynomial in 
$x$.  Let $M$ be a fixed positive integer.  Suppose there exist a constant $A\geq 1$ and a function $G(x)$ of at most 
polynomial growth such that for $1\leq m\leq 2M$ and all sufficiently large $x$, the derivatives $\phi^{(m)}(x)$ exist and 
satisfy
\begin{equation}
\label{cond1}
\left| \frac{\phi^{(m)}(x)}{m!}\right|\leq G(x)\left(\frac{A}{x}\right)^m.
\end{equation}
Assume there exist $\eta\in (0,1)$ and a constant $B$ such that for all sufficiently large $x$, 
\begin{equation}
\label{add}
G(y)\leq BG(x)\quad {\rm whenever}\ |y-x|\leq \eta x.
\end{equation}
Let $\Omega$ be an unbounded subset of $[0,\infty)$ and let $U_x,\ x\in \Omega,$ be a family of 
nonnegative random variables.  Assume there exist a constant $\delta>0$ and a function $g(s)$
 such that for all $x\in \Omega$, the mgf of $U_x$ 
exists in the interval $(-\delta, \delta)$ and satisfies
\begin{equation}
\label{mgf}
E e^{sU_x}=e^{xg(s)},\quad s\in (-\delta,\, \delta). 
\end{equation}
Assume $g'(0)=1$ in addition.  Then
\begin{equation}
\label{eqn2}
E\phi(U_x)=\phi(x)+\sum_{n=2}^{2M-2} \sum_{k=\lfloor (n+1)/2\rfloor +1}^n c_{kn} x^{n-k+1} 
\frac{\phi^{(n)}(x)}{n!}+O\left(G(x)x^{-M}\right),
\end{equation}
as $x$ tends to $\infty$, where $c_{kn}$ are constants that depend only on the function $g(s)$, and are determined by
$$E(U_x-x)^n=\sum_{k=0}^n c_{kn} x^{n-k+1},\quad n\geq 2.$$
\end{theorem}

Evidently, Theorem \ref{thm1} is the special case $g(s)=e^s-1$, except for the assumption (\ref{add}) on $G(x)$ which replaces 
(\ref{origin}).  This new assumption does not appear very restrictive as we shall see from the examples later in this section; 
however it does make the proof of Theorem \ref{main} more straightforward.  We also relax the assumption on $\phi(x)$ 
slightly by requiring only $2M$ derivatives. 

It should be emphasized that the function $g(s)$ in (\ref{mgf}) does not depend on $x$.  Also note that $g(s)$ is analytic 
in $s\in (-\delta, \delta)$ given the existence of the mgf.  Aside from the Poisson, examples of 
distribution families that satisfy (\ref{mgf}) include the binomial, negative binomial, and gamma families.  In general, 
suppose $Y_1, Y_2,\ldots$ is a sequence of independent and identically distributed (i.i.d.), nonnegative, nondegenerate 
random variables whose mgf exists in a neighborhood of zero.  Then the family of random variables $\{\sum_{k=1}^n Y_k,\ 
n=1,2,\ldots\}$ ($n$ is known as a convolution parameter) has mgf  
$$E\exp\left(s\sum_{k=1}^n Y_k\right) = \exp\left(n\log \left(E e^{sY_1}\right)\right),\quad n=1,2,\ldots,$$
which is of the form (\ref{mgf}) with $g(s)=(EY_1)^{-1}\log \left(E e^{sY_1}\right)$, if we index the family by its mean 
$x=nEY_1$.  This shows that Theorem \ref{main} is potentially applicable to a wide range of problems.

{\bf Example 1}  For a fixed $p\in (0,1)$, consider the family $U_x,\ x=p, 2p, 3p,\ldots$, where $U_x$ has the binomial 
distribution ${\rm Bi}(n, p),\ n=x/p$.  The first few central moments of $U_x$ are given by ($q=1-p$)
$$(\mu_2, \mu_3, \mu_4)=(qx,\ (q-p)qx,\ 3q^2x^2+q(1-6pq)x).$$
Since ${\rm Bi}(n,p)$ is a sum of $n$ i.i.d.\ Bernoulli$(p)$ random variables, (\ref{mgf}) is satisfied.

\begin{itemize}
\item
Given $r$ (real) and $a>0$, let $\phi(x)=G(x)=(x+a)^{-r}$.  It is easy to verify (\ref{cond1}) and (\ref{add}); thus we have
\begin{align}
\nonumber
\sum_{k=0}^n \binom{n}{k} p^k q^{n-k} (k+a)^{-r} = &(np+a)^{-r}+ \frac{qr(r+1)}{2} (np)(np+a)^{-r-2}\\
\nonumber
 &- \frac{(q-p)qr(r+1)(r+2)}{6} (np)(np+a)^{-r-3}\\
\label{item1}
 &+ \frac{q^2 r(r+1)(r+2)(r+3)}{8} (np)^2(np+a)^{-r-4}\\
\nonumber
 & +O(n^{-r-3}),
\end{align}
as $n\rightarrow\infty$. 

\item
If we let $\phi(x)=x^{-r},\ x\geq 1$ and $\phi(x)=0,\ x<1$, then we obtain an asymptotic expansion for
\begin{equation}
\label{invmom}
\sum_{k=1}^n \binom{n}{k} p^k q^{n-k} k^{-r}
\end{equation}
by simply substituting $a=0$ in the right hand side of (\ref{item1}).  When $r$ is a positive integer, (\ref{invmom}) is 
sometimes known as the $r$th inverse moment of the binomial.  The problem of inverse moments has a long history in statistics 
(see, for example, Stephan \cite{S}, Grab and Savage \cite{GS}, and David and Johnson \cite{DJ}).  More recently, expansions for (\ref{invmom}) have been considered by Marciniak and Wesolowski \cite{MW} (see also Rempala \cite{R}) for $r=1$, and by 
\v{Z}nidari\v{c} \cite{Z} for general $r$.  \v{Z}nidari\v{c} \cite{Z} also gives a brief historical account with many references. 

A special case corresponding to $r=-1/2$ is 
$$\sum_{k=0}^n \binom{n}{k} p^k q^{n-k} \sqrt{k}= \sqrt{np}\left[ 1-\frac{q/8}{np}+\frac{(q-p)q/16-15q^2/128}{(np)^2} + 
O\left(\frac{1}{n^3}\right)\right],$$
which is the binomial analog of (\ref{sqrt}) considered by Ramanujan (\cite{B}).

\item
Let $\phi(x)=\log(x+\beta)$ for a fixed $\beta>0$ and let $G(x)\equiv 1$.  We have (as $n\rightarrow\infty$)
\begin{align}
\nonumber
\sum_{k=0}^n \binom{n}{k} p^k q^{n-k} \log(k+\beta) =& \log(np+\beta)- \frac{npq}{2} (np+\beta)^{-2}\\
\label{kr}
 &+ \frac{(q-p)q}{3} (np)(np+\beta)^{-3}- \frac{3q^2 }{4} (np)^2(np+\beta)^{-4} \\
\nonumber
&+O(n^{-3}).
\end{align}
\end{itemize}

The problem of approximating the left hand side of (\ref{kr}) appears in Krichevskiy \cite{Kr} in an 
information theoretic context; see also Jacquet and Szpankowski \cite{JS1, JS2}, who give an alternative 
derivation of (\ref{kr}) using the method of analytic poissonization and depoissonization.  Flajolet 
\cite{F} also considers similar problems using singularity analysis.

{\bf Example 2} For a fixed $p\in (0, 1)$, consider the negative binomial family ${\rm NB}(n, p)$ whose probability mass 
function is $f(k; n,p)=\binom{n+k-1}{k} p^n q^k,\ k=0, 1,\ldots,$ where $q=1-p$.  The mean is $nq/p$ and the first few central 
moments are
$$(\mu_2, \mu_3, \mu_4)=\left(\frac{nq}{p^2},\ \frac{nq(1+q)}{p^3},\ 
\left[3+\frac{6}{n}+\frac{p^2}{nq}\right]\frac{(nq)^2}{p^4}\right).$$
Similar to the binomial case, as ${\rm NB}(n, p)$ is a sum of $n$ i.i.d.\ ${\rm geometric}(p)$ random variables, (\ref{mgf})
is satisfied and Theorem 2 is applicable for an appropriate $\phi(x)$. 

\begin{itemize}
\item
Take $\phi(x)=(x+a)^{-r},\ a>0$.  We have, as $n\rightarrow\infty$,
\begin{align}
\nonumber
\sum_{k=0}^\infty \binom{n+k-1}{k} p^n q^k (k+a)^{-r} = &\left(\frac{nq}{p}+a\right)^{-r}+ \frac{r(r+1)}{2p} 
\left(\frac{nq}{p}\right)\left(\frac{nq}{p}+a\right)^{-r-2}\\
\nonumber
 &- \frac{(1+q)r(r+1)(r+2)}{6p^2} \left(\frac{nq}{p}\right)\left(\frac{nq}{p}+a\right)^{-r-3}\\
 \label{nbmom}
 &+ \frac{r(r+1)(r+2)(r+3)}{8p^2} \left(\frac{nq}{p}\right)^2\left(\frac{nq}{p}+a\right)^{-r-4}\\
 \nonumber
 & +O(n^{-r-3}).
\end{align}

\item
As in the binomial case, we obtain an asymptotic expansion for 
\begin{equation}
\label{nbinv}
\sum_{k=1}^\infty \binom{n+k-1}{k} p^n q^k k^{-r}
\end{equation}
for real $r$ by substituting $a=0$ in the right hand side of (\ref{nbmom}).  Expansions for (\ref{nbinv}) have been considered by Marciniak and Wesolowski \cite{MW} and Rempala \cite{R} for the special case $r=1$, and by Wuyungaowa and Wang \cite{WW} for integer $r\geq 0$. 

\item
Let $\phi(x)=\log(x+\beta)$ for a fixed $\beta>0$.  We have 
\begin{align*}
\sum_{k=0}^\infty \binom{n+k-1}{k} p^n q^k \log(k+\beta) = &\log\left(\frac{nq}{p}+\beta\right)-\frac{1}{2p} 
\left(\frac{nq}{p}\right)\left(\frac{nq}{p}+\beta\right)^{-2}\\
 &+ \frac{(1+q)}{3p^2} \left(\frac{nq}{p}\right)\left(\frac{nq}{p}+\beta\right)^{-3}\\
&- \frac{3}{4p^2} \left(\frac{nq}{p}\right)^2\left(\frac{nq}{p}+\beta\right)^{-4} +O(n^{-3}).
\end{align*}

\end{itemize}

{\bf Example 3} Consider the gamma family ${\rm Gam}(x, 1),\ x>0,$ whose density function is $f(u; x)=u^{x-1} 
e^{-u}/\Gamma(x),\ u>0$.  The mean is $x$ and the first few central moments are
$$(\mu_2, \mu_3, \mu_4)=(x,\ 2x,\ 3x^2+6x).$$
The moment generating function is $(1-s)^{-x},\ s<1$, which is of the form (\ref{mgf}) with $g(s)=-\log(1-s)$. 

Take $\phi(x)=G(x)=x\log(x)$.  We have
$$\frac{1}{\Gamma(x)}\int_0^\infty u\log(u) u^{x-1} e^{-u}\, {\rm d}u = x\log(x) 
+\frac{1}{2}-\frac{1}{12x}+O\left(\frac{\log(x)}{x^2}\right),$$
as $x\rightarrow\infty$.  Noting $\Gamma(x+1)=x\Gamma(x)$, we may write
\begin{equation}
\label{digam}
\frac{\Gamma'(x+1)}{\Gamma(x+1)}=\log(x)+\frac{1}{2x}-\frac{1}{12x^2}+O\left(\frac{\log(x)}{x^3}\right),
\end{equation}
which is a familiar asymptotic formula for the digamma function (\cite{AS}, p.\ 259).  By expanding for one more term we can 
replace $O(x^{-3}\log(x))$ by $O(x^{-3})$ in (\ref{digam}).  A full asymptotic expansion can be recovered by applying 
(\ref{eqn2}) and using the following recursion between the central moments of ${\rm Gam}(x, 1)$ (see \cite{W}):
$$\mu_k=(k-1)(\mu_{k-1}+x\mu_{k-2}),\quad k\geq 2.$$

\section{Proof of Theorem \ref{main}}
Our proof follows Berndt \cite{B}.  In the setting of Theorem \ref{main} we have 

\begin{lemma}
\label{lem}
Let $h(u),\ u\in [0,\infty),$ be a Borel measurable function that can be bounded in absolute value by a polynomial.  Then for 
a fixed $t\in (0,1)$, both $EI(U_x<tx)h(U_x)$ and $EI(U_x>x/t) h(U_x)$ tend to 0 exponentially fast as $x$ tends to $\infty$, 
where $I(\cdot)$ is the indicator function.
\end{lemma}
{\bf Proof.}  Observe that $xg(s)$, the cumulant generating function of $U_x$, is an analytic function of $s$ (real) in a 
neighborhood of zero.  Because $g(0)=0,\ g'(0)=1$ and $t\in (0,1)$, we may choose $r, \epsilon>0$ small enough such that both 
$g(r)<r/t$ and $g(r+\epsilon)<r/t$.  Since $|h(u)|$ is bounded by a polynomial, there exists a constant $D$ such that 
$|h(u)|<e^{\epsilon u}+D$ for all $u\in [0,\infty)$.  We have 
\begin{align*}
|EI(U_x>x/t) h(U_x)| &\leq E(e^{\epsilon U_x} +D)e^{r(U_x-x/t)}\\
                   &= e^{x[g(r+\epsilon)-r/t]}+De^{x[g(r)-r/t]},
\end{align*}
which tends to zero exponentially as $x\rightarrow\infty$.  The proof for 
$EI(U_x<tx)h(U_x)$ is similar and hence omitted.

{\bf Proof of Theorem \ref{main}.}  Throughout we assume that $x$ is sufficiently large.  Define intervals $I_1=[0,\, 
(1-\eta)x),\ I_2=[(1-\eta)x,\, (1+\eta) x)$ and $I_3=[(1+\eta)x,\, \infty)$, where $\eta$ is as specified in (\ref{add}). 

By Lemma \ref{lem}, both
\begin{equation}
\label{phi1}
EI(U_x\in I_1)\phi(U_x)
\end{equation}
and 
\begin{equation}
\label{phi3}
EI(U_x\in I_3)\phi(U_x)
\end{equation}
tend to zero exponentially as $x\rightarrow\infty$.

Consider the Taylor polynomial 
$$\psi(y)=\sum_{k=0}^{2M-1} \frac{\phi^{(k)}(x)}{k!} (y-x)^k.$$
Since for any $y\in I_1$, 
$$|\psi(y)|\leq \sum_{k=0}^{2M-1} \left|\frac{\phi^{(k)}(x)}{k!}\right| x^k \equiv q(x),$$
we have
$$|EI(U_x\in I_1) \psi(U_x)|\leq q(x)EI(U_x\in I_1).$$
From (\ref{cond1}) it follows that $q(x)$ has at most polynomial growth as $x\rightarrow\infty$; by Lemma \ref{lem} we know 
that 
\begin{equation}
\label{psi1}
EI(U_x\in I_1)\psi(U_x)
\end{equation}
tends to zero exponentially as $x\rightarrow\infty$.

Similarly, for any $y\in I_3$, we have
$$|\psi(y)|\leq \sum_{k=0}^{2M-1} \left|\frac{\phi^{(k)}(x)}{k!}\right| y^k,$$
and hence
$$|EI(U_x\in I_3) \psi(U_x)|\leq \sum_{k=0}^{2M-1} \left|\frac{\phi^{(k)}(x)}{k!}\right|
EI(U_x\in I_3)U_x^k.$$
By Lemma \ref{lem}, each of $EI(U_x\in I_3)U_x^k,\ k\leq 2M-1,$ tends to zero exponentially as $x\rightarrow\infty$.   By 
(\ref{cond1}), each of $\phi^{(k)}(x)$ has at most polynomial growth as $x\rightarrow\infty$.  Overall 
\begin{equation}
\label{psi3}
EI(U_x\in I_3) \psi(U_x)
\end{equation}
tends to zero exponentially as $x\rightarrow\infty$.

For any $y\in I_2$, there exists some point $\zeta$ between $x$ and $y$ such that
\begin{align*}
|\phi(y)-\psi(y)| &=\left|\frac{\phi^{(2M)}(\zeta)}{(2M)!}\right| |y-x|^{2M}\\
                  &\leq G(\zeta) \left(\frac{A}{(1-\eta)x}\right)^{2M} |y-x|^{2M}\\
                  &\leq BG(x) \left(\frac{A}{(1-\eta)x}\right)^{2M} |y-x|^{2M},
\end{align*}
where (\ref{cond1}) and (\ref{add}) are used in the inequalities.  Letting $C=B[A/(1-\eta)]^{2M}$, we have
$$|EI(U_x\in I_2)[\phi(U_x)-\psi(U_x)]| \leq C \frac{G(x)}{x^{2M}} E |U_x-x|^{2M}.$$

We now consider the $n$th central moment of $U_x$, $\mu_n=E(U_x-x)^n$, as a function of $x$.  (Note that the 
mean of $U_x$ is $x$ as $EU_x=xg'(0)=x$.)  Expand $xg(s)$ around $s=0$ to get
$$xg(s)=\sum_{j=1}^\infty \frac{xg^{(j)}(0) s^j}{j!}.$$
Note that the coefficient $xg^{(j)}(0)$ is the $j$th cumulant of $U_x$, and, according to the well-known relation between 
central moments and cumulants (see \cite{JKK} or \cite{W}, for example)
\begin{equation}
\label{gmu}
\mu_n=\sum_{j=0}^{n-2} \binom{n-1}{j} \mu_j xg^{(n-j)}(0),\quad n\geq 2,
\end{equation}
with $\mu_0=1$ and $\mu_1=0$.  Based on (\ref{gmu}), it is easy to show by induction that $\mu_n$ is a polynomial in $x$ of 
degree at most $\lfloor n/2\rfloor$, its coefficients depending only on the function $g(s)$.  Hence, for large $x$ we have 
$E(U_x-x)^{2M}=O(x^M)$ and 
$$EI(U_x\in I_2)[\phi(U_x)-\psi(U_x)]=O(G(x)x^{-M}).$$
Combined with the exponentially small items (\ref{phi1}), (\ref{phi3}), (\ref{psi1}) and (\ref{psi3}), this gives
$$E[\phi(U_x)-\psi(U_x)]=O(G(x)x^{-M}).$$
It remains to calculate $E\psi(U_x)$.  We have, by the definition of $c_{kn}$, 
\begin{align*}
E\psi(U_x) &=\sum_{n=0}^{2M-1} E(U_x-x)^n \frac{\phi^{(n)}(x)}{n!}\\
           &=\phi(x)+\sum_{n=2}^{2M-1} \sum_{k=0}^n c_{kn} x^{n-k+1} \frac{\phi^{(n)}(x)}{n!}\\
           &=\phi(x)+\sum_{n=2}^{2M-1} \sum_{k=\lfloor (n+1)/2\rfloor+1}^n c_{kn} x^{n-k+1} \frac{\phi^{(n)}(x)}{n!}\\
           &=\phi(x)+\sum_{n=2}^{2M-2} \sum_{k=\lfloor (n+1)/2\rfloor+1}^n c_{kn} x^{n-k+1} 
\frac{\phi^{(n)}(x)}{n!}+O(G(x)x^{-M}).
\end{align*}
Note that the inner sum over $k$ is curtailed because the degree of $\mu_n$ is at most $\lfloor n/2\rfloor$, 
i.e., $c_{kn}=0$ if $k\leq \lfloor (n+1)/2 \rfloor$.  As a consequence of (\ref{cond1}), the term corresponding to $n=2M-1$ in 
the outer sum is written as $O(G(x)x^{-M})$ in the last equality.  The proof of (\ref{eqn2}) is now complete. 

\section*{Acknowledgments}
The author would like to thank an anonymous reviewer for his/her valuable comments.

\end{document}